\documentclass[10pt,twoside]{article}
\usepackage{graphicx}
\usepackage{amsmath}
\usepackage{Latex-document}

\def\Q{{\bf Q}}
\def\Z{{\bf Z}}
\def\C{{\bf  C}}
\def\H{{\bf H}}
\def\P{{\bf P}}
\def\R{{\bf R}}
\def\O{{\mathcal O}}

\title{\bf  Hilbert Schemes of Points on Surfaces\vskip 6mm}

\author{L. G\"ottsche\vspace*{-0.5cm}\thanks{Abdus Salam
International Centre for Theoretical Physics, Strada Costiera 11,
34014 Trieste, Italy. E-mail: gottsche@ictp.trieste.it}}

\date{\vspace{-8mm}}

\begin{document}

\maketitle

\thispagestyle{first} \setcounter{page}{483}

\begin{abstract}\vskip 3mm
The Hilbert scheme $S^{[n]}$ of points on an algebraic surface $S$
is a simple example of a moduli space and also a nice (crepant) resolution of
 singularities of the symmetric power $S^{(n)}$.
For many phenomena expected for moduli spaces and nice resolutions of
singular varieties it is a model case. Hilbert schemes of points
have connections to several fields of mathematics, including mo\-du\-li spaces
of sheaves, Donaldson invariants, enumerative geometry of curves, infinite
dimensional Lie algebras and vertex algebras and also to theoretical physics.
This talk will try to  give an overview over these
connections.

\vskip 4.5mm


\end{abstract}

\vskip 12mm

\setcounter{section}{-1}

\section{Introduction} \label{section 1}\setzero
\vskip-5mm \hspace{5mm}

The Hilbert scheme  $S^{[n]}$ of points on a complex projective algebraic
surface $S$ is a a parameter variety for finite subschemes of length $n$
on $S$. It is a nice (crepant) resolution of singularities of the
$n$-fold symmetric power $S^{(n)}$ of $S$. If $S$ is a K3 surface or
an abelian surface, then $S^{[n]}$ is a compact, holomorphic symplectic (thus hyperk\"ahler)
manifold. Thus $S^{[n]}$ is at the same time a basic example of a moduli
space and an example of a nice resolution of singularities of a singular
variety. There are a number of conjectures and general phenomena, many of
which originating from theoretical physics, both about moduli spaces
for objects on surfaces and about nice resolutions of singularities.
In all of these the Hilbert scheme of points can be viewed as a model case
and sometimes as the main motivating example.
Hilbert schemes of points on a surface have connections to many
topics in mathematics, including moduli spaces of sheaves and vector bundles,
Donaldson invariants, Gromov-Witten invariants and enumerative geometry of curves,
infinite dimensional Lie algebras and vertex algebras, noncommutative geometry
and also theoretical physics.

It is usually best to look at the Hilbert schemes $S^{[n]}$ for all $n$ at
the same time, and to study their invariants in terms of generating functions,
because new structures emerge this way.
For Euler numbers, Betti numbers and conjecturally for the elliptic
genus these generating functions will be modular forms and Jacobi forms.
This fits into general conjectures from physics about invariants of moduli
spaces. Also the cohomology rings of the $S^{[n]}$  for different $n$ are closely
tied together. The direct sum over $n$ of all the cohomologies is a representation
for the Heisenberg algebra modeled on the cohomology of $S$, and the
cohomology rings of the $S^{[n]}$ can be described in terms of vertex operators.
In the case that the canonical divisor of the surface $S$ is trivial, this
leads to an elementary description of the cohomology rings of the $S^{[n]}$,
which coincides with the orbifold cohomology ring of the
symmetric power, giving a nontrivial check of a conjecture
relating the cohomology ring of a nice resolution of an orbifold to the
recently defined orbifold cohomology ring.

The Hilbert schemes $S^{[n]}$ are closely related to other moduli spaces
of objects on $S$, including moduli of vector bundles and moduli of curves
e.g. via the Serre correspondence and the Mukai Fourier transform.
This leads to applications to the geometry and topology of these moduli
spaces, to Donaldson invariants, and also to formulas in the enumerative
geometry of curves on surfaces and Gromov-Witten invariants.
We want to explain some of these results and connections.
We will not attempt to give a
complete overview, but rather give a glimpse of some of the
more striking results.

\section{The Hilbert scheme of points} \label{section 2}\setzero
\setzero\vskip-5mm \hspace{5mm}

In this article $S$ will usually be a smooth projective surface
over the complex numbers. We will study the Hilbert scheme
$S^{[n]}=Hilb^n(S)$ of subschemes of length $n$ on $S$.
The points of $S^{[n]}$ correspond to finite subschemes
$W\subset S$ of length $n$, in particular a general point corresponds just to
a set of $n$ distinct points on $S$.
  $S^{[n]}$ is projective
and comes with a universal family
$Z_n(S)\subset S^{[n]}\times S$, consisting of the $(W,x)$ with $x\in W$.
An important role in applications of $S^{[n]}$ is played by the tautological
vector bundles $L^{[n]}:=\pi_*q^*(L)$ of rank $n$ on $S^{[n]}$.
Here $\pi:Z_n(S)\to S^{[n]}$ and $q:Z_n(S)\to S$ are the projections and
$L$ is a line bundle on $S$.

Closely related to  $S^{[n]}$ is the symmetric power
$S^{(n)}=S^n/{G}_n$, the quotient of $S^n$ by the action of the
symmetric group ${G}_n$.
The points of $S^{(n)}$ correspond to effective $0$-cycles
$\sum n_i [x_i]$, where the $x_i$ are distinct points of $S$ and the
sum of the $n_i$ is $n$.
The forgetful map
$$\rho:S^{[n]}\to S^{(n)},\quad  W\mapsto \sum_{x\in S} \hbox{len}(\O_{W,x}) [x]$$
is a morphism.
The symmetric power $S^{(n)}$ is singular, as for instance the
fix-locus of any transposition in ${G}_n$ has codimension $2$.
On the other hand by \cite{Fogarty}
$S^{[n]}$ is smooth  and connected of dimension $2n$
and $\rho:S^{[n]}\to S^{(n)}$ is a resolution of singularities.
In fact this is a particularly
nice resolution:
If $Y$ is a Gorenstein variety, i.e. the dualizing sheaf is a line bundle
$K_Y$,  a resolution $f:X\to Y$
of singularities is called {\it crepant} if it preserves the
canonical divisor, that is  $f^*K_Y=K_X$.
It is easy to see that $\rho:S^{[n]}\to S^{(n)}$ is crepant.
In the special case that $S$ is an abelian surface or a K3 surface
one can get a better result:
A complex manifold $X$ is called {\it holomorphic symplectic} if there
exists an everywhere non-degenerate holomorphic $2$-form $\phi$ on $X$.
If furthermore $\phi$ is unique up to scalar, $X$ is called irreducible
holomorphic symplectic.
A K\"ahler manifold $X$ of real dimension $4n$ is called {\it hyperk\"ahler}
if its holonomy group is $Sp(n)$. Compact complex manifolds are
holomorphic symplectic if and only of they admit a hyperk\"ahler metric.
In \cite{Bea1} it is shown that for a K3 surface $S$ the Hilbert scheme
$S^{[n]}$ is irreducible holomorphic symplectic. There also, for an
abelian surface $A$, the generalized
Kummer varieties are constructed from $A^{[n]}$. They form another series
of irreducible holomorphic symplectic manifolds.
The only other examples of compact hyperk\"ahler manifolds,
known not to be diffeomorphic to one in the above two series
are the two isolated examples of resolutions of singular moduli spaces of
sheaves on K3 and abelian surfaces in \cite{O'Grady1},\cite{O'Grady2}.

\section{Betti number, Euler numbers, elliptic genus}\setzero

\vskip-5mm \hspace{5mm}

For many questions about the Hilbert schemes
$S^{[n]}$ one should look at all $n$ at the same time. The first instance of this are the Betti numbers and
Euler numbers, for which we can find generating functions in terms of
modular forms.
Let ${\cal H}:=\bigl\{\tau\in \C\bigm| \Im(\tau)>0\bigr\}$.
A modular form of weight $k$ on $Sl(2,\Z)$ is a function
$f:{\cal H}\to \C$ s.th.
$$f\left( \frac{a\tau+b}{c\tau+b}\right)=(c\tau+d)^kf(\tau), \quad
\left(\begin{matrix}
a&b\\c&d\end{matrix}\right)\in Sl(2,\Z).$$
Furthermore, writing $q=e^{2\pi i\tau}$, we require that, in the Fourier
development  $f(\tau)=\sum_{n\in \Z } a_n q^n$, all the
the negative Fourier coefficients vanish. If also $a_0=0$, $f$ is called
a cusp form. The most well-known modular form is the discriminant
$\Delta(\tau):=q\prod_{n>0}(1-q^n)^{24},$
the unique cusp form of weight $12$. The Dirichlet eta function is
$\eta=\Delta^{1/24}$.

For a manifold $X$ we denote by
$p(X,z):=\sum_i (-1)^i b_i(X) z^i$ the Poincar\'e polynomial and
by $e(X)=p(X,1)$ the Euler number. The Betti numbers and Euler numbers of the
$S^{[n]}$ have very nice generating functions \cite{G1}:
\begin{equation}\sum_{n\ge 0} p(S^{[n]},z)t^n=\prod_{k\ge 1}\prod_{i=0}^4
(1-z^{2k-2+i}t^k)^{(-1)^{i+1}b_i(S)}.
\label{Gform}\end{equation}
In particular
$\sum_{n\ge 0} e(S^{[n]})q^{n-e(S)/24}=\eta(\tau)^{-e(S)}.$

This was first shown in \cite{ES1} in the case of the
projective plane and of Hirzebruch surfaces using a natural $\C^*$ action.
The proof in \cite{G1} uses the Weil Conjectures.
An important role in this proof as in all subsequent generalizations and
refinements is played by the following natural stratification of $S^{[n]}$ and
$S^{(n)}$ parametrized by the set $P(n)$ of partitions of $n$.
For a partition  $\alpha=(n_1,\ldots,n_r)\in P(n)$, the corresponding locally closed stratum
$S^{(n)}_\alpha$
of $S^{(n)}$ consists of the set of zero cycles
$n_1[x_1]+\ldots+n_r [x_r]$ with $x_1,\ldots,x_r$ distinct points of $S$.
We put $S^{[n]}_\alpha=\rho^{-1}(S^{(n)}_\alpha)$.
A partition  $\alpha=(n_1,\ldots,n_r)\in P(n)$ can also be written
as $\alpha=(1^{\alpha_1},\ldots,n^{\alpha_n})$, where
$\alpha_i$ is the number of occurences of $i$ in $(n_1,\ldots,n_r)$.
We put $|\alpha|=r=\sum \alpha_i$.
Then (\ref{Gform}) can be reformulated as
\begin{equation} p(S^{[n]},z)=\sum_{\alpha\in P(n)} p(S^{(\alpha_1)}\times\ldots\times
S^{(\alpha_n)},z)z^{2(n-|\alpha|)}.\label{hilpol}
\end{equation}
This result has been refined to Hodge numbers in \cite{GS}, \cite{Cheah}
and this was generalized in \cite{dCM1} to the Douady space of a complex
surface. It has been further refined to determine the motive
and the Chow groups \cite{dCM2} and the element in the Grothendieck
group of varieties  of $S^{[n]}$ \cite{Gmotive}.

Partially motivated by (\ref{Gform})
and using arguments from physics in \cite{DMVV}
a conjectural refinement to the Krichever-H\"ohn
elliptic genus is given.
We restrict our attention to the case that $K_X=0$ when
the elliptic genus is a Jacobi form.
For a complex vector bundle $E$ on a complex manifold $X$ and a variable $t$ we put
$$\Lambda_t(E):=\bigoplus_{k\ge 0} \Lambda^k( E) t^k, \
S_t(E):=\bigoplus_{k\ge 0} S^k( E) t^k.$$
For the holomorphic Euler characteristic we write
$\chi(X,\Lambda_t(E)):=\sum \chi(X,\Lambda^kE)t^k$ and similarly
for $S_t(E)$.
Then the elliptic genus is defined by
$$\phi(X,q,y):=\chi\left(X,\prod_{m\ge 1} \Lambda_{-y^{-1}q^m}T_X\otimes
\Lambda_{-y q^{m-1}}T^*_X\otimes S_{q^m}(T_X\oplus T_X^*)\right).$$
Writing
$\phi(S):=\sum_{m\ge 0,l} c(m,l)q^my^l$,
the conjecture  is
$$\sum_{N\ge 0} \phi(S^{[n]})p^N=\prod_{n>0,m\ge 0, l}\frac{1}{(1-p^nq^my^l)^{c(nm,l)}}.$$

\section{Infinite dimensional Lie algebras and the cohomology ring}\setzero

\vskip-5mm \hspace{5mm}

We saw that one gets nice generating functions in $n$ for the Betti numbers of the
$S^{[n]}$.
Now we shall see that the direct sum of all the cohomologies of the $S^{[n]}$
carries a new structure which governs the ring structures of the Hilbert
 schemes.
We only consider cohomology with rational coefficients and thus write
$H^*(X)$ for $H^*(X,\Q)$.
We write $H:=H^*(S)$; for $n>0$ let $\H_n:=H^*(S^{[n]})$.
and $\H:=\bigoplus_{n\ge 0} \H_n$. We shall see that
$\H$ is an irreducible module under a Heisenberg algebra.
This was conjectured in \cite{VW} and proven in  \cite{Nak1},\cite{Groj}.
$\H$ contains a distingished element $\bf 1\in H_0=\Q$.
We denote by $\int_S$ and $\int_{S^{[n]}}$ the evaluation on the fundamental
class of $S$ and $S^{[n]}$.
Define for $n>0$ the incidence
variety
$$Z_{l,n}:=\big\{(Z,x,W)\in S^{[l]}\times S\times S^{[l+n]}\bigm|
Z\subset W,\rho(W)-\rho(Z)=n[x]\bigr\},$$
 and use this to define operators
 $$p_{n}:H\to End(\H);\  p_{n}(\alpha)(y):=pr_{3*}(pr_2^*(\alpha)\cup pr_1^*(y)\cap
 [Z_{l,n}]).$$
Let $p_{-n}(\alpha):=(-1)^np_{-n}(\alpha)^\dagger$, where $\vphantom{A}^\dagger$
 denotes the adjoint with respect to $\int_{S^{[n]}}$, and $p_0(\alpha):=0$.
By \cite{Nak1},\cite{Groj} the $p_n(\alpha)$
fulfill the commutation relations of a Heisenberg algebra:
 \begin{equation}\label{Heis}
 [p_n(\alpha),p_m(\beta)]=
 (-1)^{n-1}n \delta_{n,-m}\Big(\int_S \alpha\cdot \beta\Big) id_{\H},\qquad
 n,m\in\Z,\ \alpha,\beta\in H.
 \end{equation}
We can interpret this as follows.
Let $H=H_+\oplus H_-$ be the decomposition into even and odd cohomology.
Put
$S^*(H):=\bigoplus_{i\ge 0} S^i(H_+)\otimes\bigoplus_{i\ge 0} \Lambda^i(H_-).$
The Fock space associated to $H$ is
$F(H):=S^*(H\otimes t^{}\Q[t])$.
Using the above theorems one readily shows that there is an isomorphism
of graded vector spaces
$F(H)\to \H$.
With this $\H$ becomes an irreducible module under the Heisenberg-Clifford algebra.

The ring structure of the $H^*(S^{[n]})$ is connected to the
Heisenberg algebra action. Given an action of a Heisenberg algebra,
a standard construction gives an action of the corresponding Virasoro algebra.
The important fact however, proven in \cite{Lehn}
is that the Virasoro algebra generators have a geometrical interpretation
tying them to the ring structure of the cohomology of the $S^{[n]}$.
Let $\delta:S\to S\times S$ be the diagonal embedding, and let
$\delta_*:H^*(S)\to H^*(S\times S)$ be the corresponding pushforward.
Let $p_{\nu} p_{n-\nu}\delta(\alpha):H^*(S)\to End(\H)$ be defined as
$p_\nu p_{n-\nu}(\beta\times \gamma):=p_\nu(\beta)
 p_{n-\nu}(\gamma)$ applied to $\delta_*(\alpha)\in H\times H$.
 For $n\ne 0$ define
 $L_{n}(\alpha):=\sum_{\nu\in \Z} p_{\nu}p_{n-\nu}\delta_*(\alpha),$
 and
 $L_0(\alpha):=\sum_{\nu>0} p_\nu p_{-\nu} \delta_*(\alpha).$
 These operators satisfy the relations of the
 Virasoro algebra:
 \begin{equation}\label{vir1}[L_n(\alpha),L_m(\beta)]=(n-m)L_{n+m}(ab)+\delta_{n,-m}
 \frac{n^3-n}{12}\Big(\int_S c_2(S) a b \Big) id_{\H}.
 \end{equation}
Let $\partial:\H\to \H$ be the operator which on each
$H^*(S^{[n]})$ is the multiplication with $c_1(\O^{[n]})$, where
$\O^{[n]}=\pi_*(Z_n(S))$ is the tautological vector bundle associated to the
trivial line bundle on $S$.
The tie given in  \cite{Lehn}  to the ring structure is:
\begin{equation}\label{vir2}
[\partial, p_n(\alpha)]=nL_n(a)+{n \choose 2} p_n(K_S a),\qquad
n\ge 0,\ \alpha\in H^*(S) .\end{equation}
In \cite{LQW1}, for each $\alpha\in H^*(S)$, classes $\alpha^{[n]}\in H^*(S^{[n]})$
are defined as generatizations of the Chern characters $ch(F^{[n]})$
of tautological bundles, which are studied in \cite{Lehn}.
The homogeneous components of the $\alpha^{[n]}$ generate
the ring $H^*(S^{[n]})$.
\cite{Lehn},\cite{LQW1} relate the multiplication by the $\alpha^{[n]}$
to the higher order commutators with $\partial$: Let
$\alpha^{[\bullet]}:\H\to \H$ be the operator which on every $H^*(S^{[n]})$ is
the multiplication with $\alpha^{[n]}$, then
\begin{equation}\label{vir3}
[\alpha^{[\bullet]},p_1(\beta)]=\exp(ad(\partial)) p_1(\alpha\beta),
\end{equation}
where for an operator $A:\H\to \H$, $ad(\partial)A=[\partial, A]$.

(\ref{vir1}),(\ref{vir2}),(\ref{vir3}) determine the cohomology rings
of the $S^{[n]}$.
In case $K_S=0$ this is used in \cite{LS1},\cite{LS2}  to give an
elementary description of the cohomology rings $H^*(S^{[n]})$ in terms of
the symmetric group, which
we will relate below to orbifold cohomology rings.

\section{Orbifolds and orbifold cohomology}\setzero

\vskip-5mm \hspace{5mm}

Let  $X$ be a compact complex manifold with an action of a finite
group $G$ and assume that for all $1\ne g\in G$
the fixlocus $X^g$
has codimension $\ge 2$. The quotient $X/G$ will usually be singular, but
the stack quotient $[X/G]$  is a smooth
orbifold. In physics \cite{DHVW1},\cite{DHVW2} the following orbifold
Euler characteristic has been introduced
$$e(X,G):=\sum_{gh=hg\in G} e(X^{g,h})=
\sum_{[g]\subset G}e(X^g/C(g)).$$
Here the first sum runs over all commuting pairs in $G$ and $X^{g,h}$
is the set of common fixpoints; the second sum runs over the conjugacy
classes $[g]$ of elements in $G$ and $C(g)$ is the centralizer of $g$.
If $Y\to X/G$ is a crepant resolution, then it was expected that
$e(X,G)=e(Y)$. As the conjugacy classes of the
symmetric group $G_n$ correspond to the partitions of $n$, one can
see \cite{HH} using formula (\ref{hilpol}) that
this is true for the resolution of $S^{(n)}$ by $S^{[n]}$, which
was an important  check for this conjecture.

Orbifold Euler numbers have been refined to
 orbifold cohomology groups \cite{Z}. We again take all cohomology with $\Q$
coefficients.
Define a rationally graded $\Q$-vector space
$$H^*_{orb}([X/G]):=\bigoplus_{[g]\subset G} H^*(X^g/C(g)).$$
The grading is defined as follows.
Assume for simplicity that all $X^g$ are connected.
 For $x\in X^g$ let $e^{2\pi i r_1},\ldots,e^{2\pi i r_k}$ be the
 eigenvalues of $g$ on $T_{X,x}$.
 Put $a(g):=\sum r_i\in \Q$ where $r_i\in [0,1)$. This
 is independent of $x$. For  $\alpha\in H^i(X^g/C(g))$ its degree in
 the $[g]$-th summand of $H^*_{orb}([X/G])$ is $i+2a(g)$.
 If $X/G$ is Gorenstein,
 then it is easy to see that $a(g)\in \Z_{\ge 0}$.
 For crepant resolutions $Y\to X/G$,
 it was conjectured that $H^*_{orb}([X/G])=H^*(Y)$ as graded
 vector spaces.
In the case of  $S^{[n]}\to S^{(n)}$ this can again be verified from
formula (\ref{hilpol}).
In  \cite{BB} it has been established for all
crepant resolutions $Y\to X/G$.

 Recently orbifold cohomology rings, i.e. a ring structure
 on the orbifold cohomology have been defined as a special case of
quantum cohomology of orbifolds \cite{CR1},\cite{AV},\cite{AGV}.
 In \cite{R1} it is conjectured for  an orbifold $X$ with a
hyperk\"ahler resolution  $Y\to X$, i.e. a crepant resolution
 such that $Y$ is hyperk\"ahler, that the orbifold cohomology ring of $X$ is
  isomorphic to  $H^*(Y)$. The most relevant case of such
 a resolution is  $S^{[n]}\to S^{(n)}$ when $K_S=0$.
This is precisely the case in which \cite{LS2} gives an elementary
description of the cohomology ring of $S^{[n]}$. In \cite{FG},\cite{U}
an elementary description of the orbifold cohomology of
a quotient $[X/G]$ by a finite group is given.
We define  $H^*(X,G):=\sum_{g\in G} H^*(X^g)$. This carries a
$G$-action by $h\cdot\alpha_g=(h_*\alpha)_{hgh^{-1}}$, and, for a suitable grading
on  $H^*(X,G)$, it follows that
the $G$ invariant part is just $H^*_{orb}([X/G])$ as a graded vector space.
In order to define the ring structure on  $H^*_{orb}([X/G])$ one therefore
defines a ring structure on $H^*(X,G)$ compatible with the $G$-action.
In \cite{LS2} the cohomology ring $H^*(S^{[n]})$ is also described as
the $G_n$ invariant part of a ring structure on $H^*(S^n,G_n)$ and one checks
that the two ring structures on $H^*(S^n,G_n)$ coincide up to an explicit sign change,
thus proving the conjecture of  \cite{R1} for $S^{[n]}$.

If $\pi:Y\to X/G$ is only a crepant resolution but not hyperk\"ahler,
  then usually
 $H^*_{orb}([X/G])$ and $ H^*(Y)$ are not isomorphic as rings.
However in \cite{R2} a precise conjecture is made relating the two:
One has to correct   $H^*_{orb}([X/G])$ by Gromov-Witten invariants
coming from  classes of rational curves $Y$ contracted
by $\pi$. In the case of the Hilbert scheme
these curve classes are the multiples of a unique class.
The conjecture was verified  for $S^{[2]}$.

\section{Moduli of vector bundles}\setzero

\vskip-5mm \hspace{5mm}

We denote by $M_S^H(r,c_1,c_2)$ the  moduli space of Gieseker $H$-semistable
 coherent sheaves
of rank $r$ on $S$ with Chern classes $c_1$, $c_2$. Here a sheaf
$\cal F$ of rank $r>0$ on $S$ is called semistable, if
$\chi({\cal G}\otimes H^n)/r'\le \chi({\cal F}\otimes H^n) /r$ for all sufficiently
large $n$ and for all subsheaves
$\cal G\subset F$ of positive rank $r'$.
As
$M_S^H(1,0,c_2)\simeq Pic^0(S)\times S^{[c_2]}$, the Hilbert scheme of
points is a special case.
We will often restrict our attention to the case of $r=2$ and
 write  $M_S^H(c_1,c_2)$.
The Hilbert schemes of points are related in several ways to the
$M_S^H(c_1,c_2)$. The most basic tie is the Serre correspondence which says
that under mild assumptions rank two vector bundles on $S$ can be constructed
as extensions of ideal sheaves of finite subschemes by line bundles.
Related to this is the dependence of the $M_S^H(c_1,c_2)$ on the
ample divisor $H$ via a system of walls and chambers. This has been studied by
a number of authors (e.g. \cite{Qin},\cite{FQ},\cite{EG}).
Assume for simplicity that $S$ is simply connected.
A class $\xi\in H^2(S,\Z)$ defines a wall of type $(c_1,c_2)$ if $\xi+c_1\in
2H^2(S,\Z)$ and $c_1^2-4c_2\le \xi^2<0$. The corresponding wall is
$W^\xi=\big\{\alpha\in H^2(S,\R)\bigm|\alpha\cdot \xi=0\big\}$.
The connected components of the complement of the walls in $H^2(X,\R)$
are called the chambers of type $(c_1,c_2)$.
If a sheaf ${\cal E}\in M_S^H(c_1,c_2)$ is unstable with respect to $L$,
then there is a wall $W^\xi$ with $H\xi<0<L\xi$ and an extension
$$0\to {\cal I}_{Z}\otimes A\to {\cal E}\to {\cal I}_{W}\times B\to 0,$$
where $A,B\in Pic(S)$ with $A-B=\xi$ and ${\cal I}_{Z}$, ${\cal I}_{W}$
are the ideal sheaves of zero dimensional schemes on $S$.
It follows that
$M_S^H(c_1,c_2)$ depends only on the chamber of $H$ and the set theoretic
change under wallcrossing is given in terms of Hilbert schemes of points on $S$.
In the case e.g. of rational surfaces and K3-surfaces, the change can
be described as
an explicit sequence  of blow ups along $\P_k$ bundles over products
$S^{[n]}\times S^{[m]}$ followed by blow downs in another direction
\cite{FQ},\cite{EG}.
The change of the Betti and Hodge numbers under wallcrossing can be explicitly determined and this can
be used e.g. to determine the Hodge numbers of $M_S^H(c_1,c_2)$ for
rational surfaces. For suitable choices of $H$ one can find the generating
functions in terms of modular forms and Jacobi forms \cite{Gtheta}.

The appearance of modular forms is in accord with the $S$-duality
conjectures \cite{VW} from theoretical physics, which predict
that under suitable assumptions the generating functions for the
the Euler numbers of moduli spaces of sheaves on surfaces should
be given by modular forms. One of the motivating examples
for this conjecture is the case that $S$ is a K3-surface.
In this case the  conjecture is that, if  $M_S^H(c_1,c_2)$ is smooth,
then it has the same Betti numbers as the Hilbert scheme of points
on $S$ of the same dimension. Assuming this, the formula (\ref{Gform})
for the Hilbert schemes of points
implies that the generating function for the Euler numbers is a modular
form. If $c_1$ is primitive this was shown in \cite{GH}.
The result was shown in general
for
$M_S^H(r,c_1,c_2)$ with $r>0$ in \cite{YK31},\cite{YK32},
by relating the Hilbert scheme and the
moduli space via birational correspondences and deformations.
One concludes that $M_S^H(r,c_1,c_2)$ has the same Betti numbers
as the Hilbert scheme of points of the same dimension, as both
spaces are holomorphic symplectic \cite{Muk} and birational
manifolds with trivial canonical class have the same Betti numbers
\cite{Ba}.
 Similar results are shown in \cite{YA} for abelian surfaces.
 Other motivating examples for the S-duality conjecture were the
 case of
$\P_2$ \cite{YP2} and the blowup formula relating the
generating function for the Euler numbers of the moduli spaces
of rank 2 sheaves on a surface $S$ to that on the
blowup of $S$ in a point, which has since been established
(\cite{LQ1},\cite{LQ2}, see also \cite{Gtheta}).

The moduli spaces $M_S^H(c_1,c_2)$ can be used to compute the Donaldson
invariants of $S$. In case $p_g=0$ these depend on a metric, corresponding to the dependence of
$M_S^H(c_1,c_2)$ on $H$. For rational surfaces one can use the above
description of the wallcrossing for the $M_S^H(c_1,c_2)$ to determine the change
of the Donaldson invariants in terms of Chern numbers of  generalizations of the
tautological sheaves $L^{[n]}$ on products $S^{[n]}\times S^{[m]}$
of Hilbert schemes of points \cite{EG},\cite{FQ}.
The leading terms of these expressions can be
explicitly evaluated. The wallcrossing of Donaldson invariants has also been
studied in gauge theory (e.g.\cite{K},\cite{KM}).
There a conjecture about the structure
of the wallcrossing  formulas is made. Assuming this conjecture one can determine
the generating functions for the wallcrossing in terms of modular forms
\cite{Gmod},\cite{GZ}.

\section{Enumerative geometry of curves}\setzero

\vskip-5mm \hspace{5mm}

Now we want to see some striking relations between the Hilbert schemes
$S^{[n]}$ and the enumerative geometry of curves on $S$.
First let $S$ be a K3 surface and  $L$ a primitive line bundle on $S$.
Then $L^2=2g-2$, where the linear
system $|L|$ has dimension $g$ and a smooth curve in $|L|$ has geometric genus
$g$. As a node imposes one linear condition, one
expects a finite number of rational
curves (i.e. curves of geometric genus $0$) in $|L|$.
 Partially based on arguments  from physics, a formula
is given in \cite{YZ} for the number of rational curves in $|L|$ and in
\cite{Beak3} this made mathematically precise.
Writing $n_g$ for the number of rational curves in $|L|$ with $L^2=2g-2$
(counted with suitable multiplicities), the formula is
\begin{equation}
\sum_{g\ge 0}{n_g}q^g=\frac{q}{\Delta},\label{k3curv}
\end{equation}
where $\Delta$ is again the discriminant.
By (\ref{Gform}) this implies the surprizing fact that
$n_g$ is just the Euler number of $S^{[g]}$. In fact
the argument  relates the number of curves to $S^{[g]}$:
Let ${\mathcal C}\to |L|$ be the universal curve and let ${\mathcal J}\to |L|$
be the corresponding relative compactified Jacobian, whose fibre over the
point corresponding to a curve $C$ is the compactified Jacobian $J(C)$
\cite{AK1}. One can show that $e(J(C))=0$ unless $g(C)=0$.
It follows that
$e({\mathcal J})$ is the sum over the $e(J(C))$ for $C\in |L|$ with
$g(C)=0$. It is not difficult to show that $S^{[g]}$ and ${\mathcal J}$ are
birational. $\mathcal J$ is also smooth and hyperk\"ahler as a moduli
space of sheaves on a K3 surface \cite{Muk}. As already used in the
section on vector bundles, birational manifolds with
trivial canonical bundle have the same Betti numbers \cite{Ba}.
Thus
$\mathcal J$ and $S^{[g]}$ have  the same Euler numbers.  This shows
(\ref{k3curv}), where the multiplicity of a rational curve $C$
is $e(J(C))$. By \cite{FGS}
this multiplicity is the multiplicity of the
corresponding moduli space of stable maps, in particular it is
always positive.
In \cite{Gconj} a conjectural generalization of (\ref{k3curv}) to arbitrary surfaces
$S$ is given.

\smallskip

\noindent{\bf Conjecture 6.1}{\it
\begin{enumerate}
\item For all $\delta\ge 0$, there exists a universal polynomial
$T_\delta(x,y,z,w)$, such that
for all projective surfaces $S$ and all sufficiently ample
line bundles $L$ on $S$ the number of $\delta$-nodal curves in a general
$\delta$-dimensional linear subspace of $|L|$ is
$T_{\delta}(\chi(L),\chi(\O_S),LK_S,K_S^2)$.
\item
There are universal power series $B_1,B_2\in \Z[[q]]$ whose coefficients can
be explicitely determined, such that
$$\sum_{\delta\ge 0}
T_{\delta}(\chi(L),\chi(\O_S),LK_S,K_S^2)(DG_2)^\delta=
\frac{(DG_2/q)^{\chi(L)}B_1^{LK_S}B_2^{K_S^2}}
{(\Delta D^2 G_2/q^2)^{\chi(\O_S)/2}}.
$$
Here $D=q\frac{d}{dq}$ and $G_2=-\frac{1}{24}\frac{D\Delta}{\Delta}$.
\end{enumerate}}

The expectation that universal polynomials  should exist is
implicit in \cite{Vain},\cite{KP1} where the $T_\delta$ are
determined for $\delta\le 8$. In \cite{Gconj} also another tie of
the conjecture to the Hilbert scheme of points is given:
conjecturally the numbers
$T_{\delta}(\chi(L),\chi(\O_S)$,\\$LK_S,K_S^2)$ are suitable
intersection numbers on the Hilbert scheme $S^{[3\delta]}$ of
$3\delta$ points of $S$. If $S$ is a K3 surface or an abelian
surface, then the conjecture predicts that the generating function
can be written in terms of modular forms. In this case a modified
version of Conjecture 6.1 was proven for primitive line bundles in
\cite{Bryk3} and \cite{Bryab}, replacing the numbers of
$\delta$-nodal curves with the corresponding modified
Gromov-Witten invariants. In \cite{Liu} a proof of the conjecture
is published.

\label{lastpage}

\end{document}